\documentclass[11pt]{article}

\usepackage{amssymb}
\usepackage{amsmath,amssymb,epsf,wrapfig,multicol,epic,eepic,trig,mathrsfs,dsfont,color,graphicx}
 \usepackage[abs]{overpic}
\def\qed{\hfill$\square$}
\newtheorem{thm}{Theorem}

\newtheorem{prop}[thm]{Proposition}

\def\hsymbu#1{\smash{\lower1.7ex\hbox{\huge$#1$}}}

\def\rmoveio#1#2{%#2=2 oriented #2=1 non-oriented
\setlength{\unitlength}{#1}
\begin{picture}(50,30)
\put(5,0){\line(0,1){30}}

{\allinethickness{.8pt}
\put(10,15){\vector(1,0){15}}
\put(25,15){\vector(-1,0){15}}}

\qbezier(30,0)(30,20)(45,20)
\qbezier(45,20)(50,20)(50,15)
\qbezier(50,15)(50,10)(45,10)
\qbezier(45,10)(40,10)(36,14)
\qbezier(33,17)(30,25)(30,30)

\ifnum#2=2
\put(3,28){\path(0,0)(2,2)(4,0)}
\put(28,28){\path(0,0)(2,2)(4,0)}
\put(38,15){\makebox{${\Huge c_{1}}$}}
\fi

\end{picture}
}

\def\rmoveiio#1#2{%#2=2 oriented #2=1 non-oriented #2=3 cohearent oriented
\setlength{\unitlength}{#1}
\begin{picture}(60,40)
\put(5,0){\line(0,1){40}}
\put(15,0){\line(0,1){40}}

{\allinethickness{.8pt}
\put(20,20){\vector(1,0){15}}
\put(35,20){\vector(-1,0){15}}}

\qbezier(40,0)(42,1)(47,5)
\qbezier(50,8)(68,20)(50,32)
\qbezier(47,35)(42,39)(40,40)

\qbezier(60,0)(20,20)(60,40)

\ifnum#2=2
\put(2,37){\path(0,0)(3,3)(6,0)}
\put(12,37){\path(0,0)(3,3)(6,0)}
\put(40,37){\path(0,0)(0,3)(3,3)}
\put(60,37){\path(0,0)(0,3)(-3,3)}
\put(47,25){\makebox{${\Huge c_{1}}$}}
\put(47,13){\makebox{${\Huge c_{2}}$}}
\fi

\ifnum#2=3
\put(2,2){\path(0,0)(3,-3)(6,0)}
\put(12,37){\path(0,0)(3,3)(6,0)}
\put(40,3){\path(0,0)(0,-3)(3,-3)}
\put(60,37){\path(0,0)(0,3)(-3,3)}
\put(47,25){\makebox{${\Huge c_{1}}$}}
\put(47,13){\makebox{${\Huge c_{2}}$}}
\fi

\end{picture}
}
\def\rmoveiiio#1#2{%#2=2 oriented #2=1 non-oriented #2=3 cohearent oriented
\setlength{\unitlength}{#1}
\begin{picture}(75,30)
\put(0,0){\line(1,1){15}}
\qbezier(15,15)(20,20)(20,30)

\put(10,0){\line(-1,1){4}}
\qbezier(4,6)(-5,15)(5,25)
\put(5,25){\line(1,1){5}}

\qbezier(20,0)(20,10)(16,14)
\put(14,16){\line(-1,1){8}}
\put(4,26){\line(-1,1){4}}

{\allinethickness{.8pt}
\put(30,15){\vector(1,0){10}}
\put(40,15){\vector(-1,0){10}}}

\qbezier(50,0)(50,10)(55,15)
\put(55,15){\line(1,1){15}}

\put(60,0){\line(1,1){5}}
\qbezier(65,5)(75,15)(66,24)
\put(64,26){\line(-1,1){4}}

\put(70,0){\line(-1,1){4}}
\put(64,6){\line(-1,1){8}}
\qbezier(54,16)(50,20)(50,30)

\ifnum#2=2
\put(0,27){\path(0,0)(0,3)(3,3)}
\put(7,30){\path(0,0)(3,0)(3,-3)}
\put(17,27){\path(0,0)(3,3)(6,0)}

\put(10,23){\makebox{${\Huge c_{1}}$}}
\put(10,3){\makebox{${\Huge c_{2}}$}}
\put(20,13){\makebox{${\Huge c_{3}}$}}

\put(47,27){\path(0,0)(3,3)(6,0)}
\put(60,27){\path(0,0)(0,3)(3,3)}
\put(67,30){\path(0,0)(3,0)(3,-3)}

\put(70,23){\makebox{${\Huge c'_{2}}$}}
\put(70,3){\makebox{${\Huge c'_{1}}$}}
\put(60,13){\makebox{${\Huge c'_{3}}$}}
\fi

\end{picture}
}
\def\rmovevio#1#2{%#2=2 oriented #2=1 non-oriented
\setlength{\unitlength}{#1}
\begin{picture}(50,30)
\put(5,0){\line(0,1){30}}

{\allinethickness{.8pt}
\put(10,15){\vector(1,0){15}}
\put(25,15){\vector(-1,0){15}}}

\qbezier(30,0)(30,20)(45,20)
\qbezier(45,20)(50,20)(50,15)
\qbezier(50,15)(50,10)(45,10)
\qbezier(45,10)(30,10)(30,30)

\put(34,15){\circle{5}}

\ifnum#2=2
\put(3,28){\path(0,0)(2,2)(4,0)}
\put(28,28){\path(0,0)(2,2)(4,0)}
\put(38,15){\makebox{${\Huge c_{1}}$}}
\fi

\end{picture}
}

\def\rmoveviio#1#2{%#2=2 oriented #2=1 non-oriented #2=3 cohearent oriented
\setlength{\unitlength}{#1}
\begin{picture}(60,40)
\put(5,0){\line(0,1){40}}
\put(15,0){\line(0,1){40}}

{\allinethickness{.8pt}
\put(20,20){\vector(1,0){15}}
\put(35,20){\vector(-1,0){15}}}

\qbezier(40,0)(80,20)(40,40)
\qbezier(60,0)(20,20)(60,40)
\put(50,6){\circle{5}}
\put(50,34){\circle{5}}

\ifnum#2=2
\put(2,37){\path(0,0)(3,3)(6,0)}
\put(12,37){\path(0,0)(3,3)(6,0)}
\put(40,37){\path(0,0)(0,3)(3,3)}
\put(60,37){\path(0,0)(0,3)(-3,3)}
\put(47,25){\makebox{${\Huge c_{1}}$}}
\put(47,13){\makebox{${\Huge c_{2}}$}}
\fi

\ifnum#2=3
\put(2,2){\path(0,0)(3,-3)(6,0)}
\put(12,37){\path(0,0)(3,3)(6,0)}
\put(40,3){\path(0,0)(0,-3)(3,-3)}
\put(60,37){\path(0,0)(0,3)(-3,3)}
\put(47,25){\makebox{${\Huge c_{1}}$}}
\put(47,13){\makebox{${\Huge c_{2}}$}}
\fi

\end{picture}
}

\def\rmoveviiio#1#2{%#2=2 oriented #2=1 non-oriented #2=3 cohearent oriented
\setlength{\unitlength}{#1}
\begin{picture}(70,30)
\put(0,0){\line(1,1){15}}
\qbezier(15,15)(20,20)(20,30)

\put(10,0){\line(-1,1){5}}
\qbezier(5,5)(-5,15)(5,25)
\put(5,25){\line(1,1){5}}

\qbezier(20,0)(20,10)(15,15)
\put(15,15){\line(-1,1){15}}
\put(5,5){\circle{5}}
\put(15,15){\circle{5}}
\put(5,25){\circle{5}}

{\allinethickness{.8pt}
\put(30,15){\vector(1,0){10}}
\put(40,15){\vector(-1,0){10}}}

\qbezier(50,0)(50,10)(55,15)
\put(55,15){\line(1,1){15}}

\put(60,0){\line(1,1){5}}
\qbezier(65,5)(75,15)(65,25)
\put(65,25){\line(-1,1){5}}

\put(70,0){\line(-1,1){15}}
\qbezier(55,15)(50,20)(50,30)

\put(65,5){\circle{5}}
\put(55,15){\circle{5}}
\put(65,25){\circle{5}}

\ifnum#2=2
\put(0,27){\path(0,0)(0,3)(3,3)}
\put(7,30){\path(0,0)(3,0)(3,-3)}
\put(17,27){\path(0,0)(3,3)(6,0)}

\put(20,13){\makebox{${\Huge c_{1}}$}}

\put(47,27){\path(0,0)(3,3)(6,0)}
\put(60,27){\path(0,0)(0,3)(3,3)}
\put(67,30){\path(0,0)(3,0)(3,-3)}

\put(60,13){\makebox{${\Huge c'_{1}}$}}
\fi

\end{picture}
}

\def\rmovevivo#1#2{%#2=2 oriented #2=1 non-oriented #2=3 cohearent oriented
\setlength{\unitlength}{#1}
\begin{picture}(70,30)
\put(0,0){\line(1,1){15}}
\qbezier(15,15)(20,20)(20,30)

\put(10,0){\line(-1,1){5}}
\qbezier(5,5)(-5,15)(5,25)
\put(5,25){\line(1,1){5}}

\qbezier(20,0)(20,10)(16,14)
\put(14,16){\line(-1,1){14}}
\put(5,5){\circle{5}}
\put(5,25){\circle{5}}

{\allinethickness{.8pt}
\put(30,15){\vector(1,0){10}}
\put(40,15){\vector(-1,0){10}}}

\qbezier(50,0)(50,10)(55,15)
\put(55,15){\line(1,1){15}}

\put(60,0){\line(1,1){5}}
\qbezier(65,5)(75,15)(65,25)
\put(65,25){\line(-1,1){5}}

\put(70,0){\line(-1,1){14}}
\qbezier(54,16)(50,20)(50,30)

\put(65,5){\circle{5}}
\put(65,25){\circle{5}}

\ifnum#2=2
\put(0,27){\path(0,0)(0,3)(3,3)}
\put(7,30){\path(0,0)(3,0)(3,-3)}
\put(17,27){\path(0,0)(3,3)(6,0)}

\put(20,13){\makebox{${\Huge c_{1}}$}}

\put(47,27){\path(0,0)(3,3)(6,0)}
\put(60,27){\path(0,0)(0,3)(3,3)}
\put(67,30){\path(0,0)(3,0)(3,-3)}

\put(60,13){\makebox{${\Huge c'_{1}}$}}
\fi

\end{picture}
}
\def\rmoveti#1{
\setlength{\unitlength}{#1}
\begin{picture}(60,20)
\put(0,10){\line(1,0){20}}
\put(15,0){\line(0,1){20}}
\put(15,10){\circle{3}}
{%\linethickness{2pt}
\put(7,7){\line(0,1){6}}} %bar
{\allinethickness{.8pt}
\put(25,10){\vector(1,0){10}}
\put(35,10){\vector(-1,0){10}}}

\put(40,10){\line(1,0){20}}
\put(45,0){\line(0,1){20}}
\put(45,10){\circle{3}}
{%\linethickness{2pt}
\put(53,7){\line(0,1){6}}}%bar
\end{picture}
}

\def\rmovetiio#1{
\setlength{\unitlength}{#1}
\begin{picture}(40,20)
\put(5,0){\line(0,1){20}}
{%\linethickness{2pt}
\put(2.,7){\line(1,0){6}} %bar
\put(2.,13){\line(1,0){6}} }%bar
{\allinethickness{.8pt}
\put(15,10){\vector(1,0){10}}
\put(25,10){\vector(-1,0){10}}}
\put(35,0){\line(0,1){20}}
\end{picture}
}

\def\rmovetiiio#1#2{%#2=2 oriented #2=1 non-oriented #2=3 cohearent oriented
\setlength{\unitlength}{#1}
\begin{picture}(70,20)

\put(0,0){\line(1,1){20}}
\put(20,0){\line(-1,1){9}}
\put(9,11){\line(-1,1){9}}

\put(7,3){\line(-1,1){4}}
\put(13,3){\line(1,1){4}}
\put(3,13){\line(1,1){4}}
\put(17,13){\line(-1,1){4}}

{\allinethickness{.8pt}
\put(25,10){\vector(1,0){10}}
\put(35,10){\vector(-1,0){10}}}

\qbezier(40,5)(48,20)(54,11)
\qbezier(56, 9)(62,0)(70,15)

\qbezier(40,15)(48,0)(55,10)
\qbezier(55,10)(62,20)(70,5)

\put(43,10){\circle{5}}
\put(67,10){\circle{5}}

\ifnum#2=2
\put(0,17){\path(0,0)(0,3)(3,3)}
\put(17,20){\path(0,0)(3,0)(3,-3)}

\put(15,8){\makebox{${\Huge c_{1}}$}}

\put(40,12){\path(0,0)(0,3)(3,3)}
\put(67,15){\path(0,0)(3,0)(3,-3)}

\put(53,3){\makebox{${\Huge c'_{1}}$}}
\fi

\end{picture}
}

\begin{document}
%% title %%%%%%%%%%%%%%%%%%%%%%%%%%%%%%
%\begin{frontmatter}
\title
{A twisted link invariant  derived from a virtual link  invariant}
\author{Naoko Kamada %and Seiichi Kamada
\thanks{This work was supported by JSPS KAKENHI Grant Number 15K04879.}
%, 26287013.}
\\ Graduate School of Natural Sciences,  Nagoya City University\\ 
1 Yamanohata, Mizuho-cho, Mizuho-ku, Nagoya, Aichi 467-8501 Japan\\
%Department of Mathematics, Osaka City University, \\
%Suimiyoshi,  Osaka 558-8585, Japan
}

\date{}
\maketitle
\begin{abstract} 
Virtual knot theory is a generalization of knot theory which is based on  Gauss chord diagrams and link diagrams on closed oriented surfaces. 
A twisted knot is a generalization of a virtual knot, which corresponds to a link diagram on a possibly non-orientable surface. 
In  this paper, we discuss an invariant of twisted links which is obtained from the JKSS invariant of virtual links by use of double coverings. We also discuss some properties of double covering diagrams.
\end{abstract}

%\begin{keyword}
%Virtual knot theory \sep Miyazawa polynomial
%\MSC Primary 57M25 \sep Secondary 57M27.
%\end{keyword}
%\end{frontmatter}

\section{Introduction}

%%%  1 virtual Link diagram%%%%%%%%%%%%%%%%%%%%%%%%
L. H. Kauffman introduced virtual knot theory, which is a generalization of knot theory based on  Gauss chord diagrams and link diagrams in closed oriented surfaces  \cite{rkauD}. 
%Virtual links correspond to  stable equivalence classes of  links in oriented 3-manifolds which are line %bundles over closed  oriented surfaces  \cite{rkk}, \cite{rCKS}. 
Twisted knot  theory was introduced by Bourgoin. It is an extension of virtual knot theory. Twisted links correspond to stable equivalence classes of links in oriented 3-manifolds which are line bundles over (possibly non-orientable) closed surfaces  \cite{rBor}, \cite{rCKS}. 

F. Jaeger, L. H. Kauffman and H. Saleur defined an invariant of links in thickened surfaces, where the surfaces are oriented \cite{rJKS}. J. Sawollek \cite{rSaw} applied it to virtual links, which is so-called the {\it JKSS invariant}.  
In this paper, we introduce an invariant of twisted links obtained from the JKSS invariant by use of double coverings.
We also discuss some properties of double coverings of twisted link diagrams.

\begin{figure}[h]
\centerline{
\includegraphics[width=6cm]{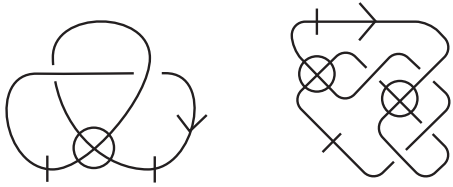}}
\caption{Examples of twisted link diagrams}\label{fig:extwtdiag}
\end{figure}

A {\it virtual link diagram\/} is a generically immersed loops whose double points have information of positive, negative or virtual crossing.  A {\it virtual crossing\/} is an encircled double point without over-under information. A {\it twisted link diagram\/} is a virtual link diagram, possibly with {\it bars\/} on arcs. Examples of twisted link  diagrams are depicted in Figure~\ref{fig:extwtdiag}.

\begin{figure}[!h]
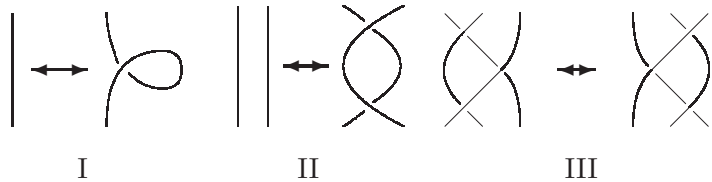

\begin{center}
\begin{tabular}{ccc}
\rmoveio{.5mm}{1}&\rmoveiio{.4mm}{1}&\rmoveiiio{.5mm}{1}\\
I&II&III\\
\end{tabular}
\caption{Reisdemeister moves}\label{fig:movesR}
\end{center}
\end{figure}

\begin{figure}[!h]
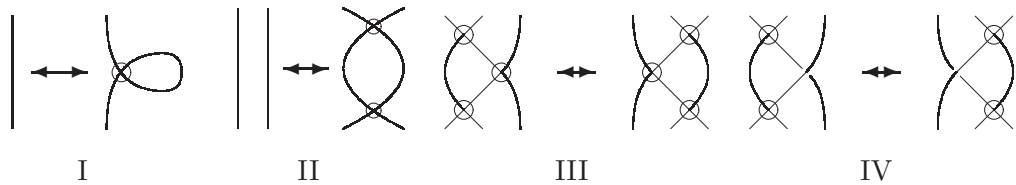

\begin{center}
\begin{tabular}{cccc}
\rmovevio{.5mm}{1}&\rmoveviio{.4mm}{1}&\rmoveviiio{.5mm}{1}&\rmovevivo{.5mm}{1}\\
I&II&III&IV\\
\end{tabular}
\caption{Virtual Reidemeister moves}\label{fig:movesV}
\end{center}
\end{figure}

\begin{figure}[!h]
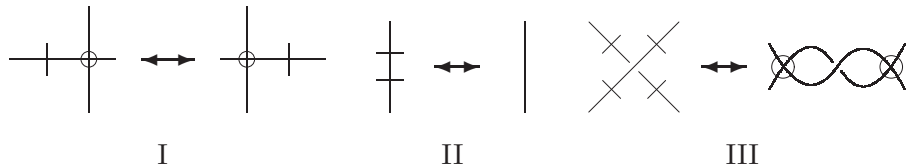

\begin{center}
\begin{tabular}{ccc}
\rmoveti{.7mm}&\rmovetiio{.6mm}&\rmovetiiio{.6mm}{1}\\
I&II&III\\
\end{tabular}
\caption{Twisted Reidemeister moves}\label{fig:movesT}
\end{center}
\end{figure}

A {\it twisted link\/}  (resp. a {\it virtual link}) is an equivalence class of a twisted link diagram (resp. virtual link diagram) under Reidemeister moves, virtual Reidemeister moves and twisted Reidemeister moves (resp. Reidemeister moves and virtual Reidemeister moves) depicted in Figures~\ref{fig:movesR}, \ref{fig:movesV} and \ref{fig:movesT}.

\section{The JKSS invariant}
\begin{figure}[h]
\centerline{
\includegraphics[width=2.5cm]{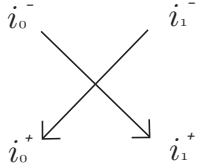}}
\caption{Labels of four edges}\label{fig:labeledge}
\end{figure}

We recall the definition of the JKSS invariant of a virtual link.
Let $D$ be a virtual link diagram with $n$ real crossings.
 Let $c_1, \dots, c_n$ be the real crossings of $D$. 
We define a $2n\times 2n$ matrix,  $M=\mathrm{diag}(M_1,\dots, M_n)$, where 
$M_i=M_+$ (or $M_-$) if $c_i$ is positive (or negative) crossing. Here $M_+$ and $M_-$
are $2\times 2$ matrices:
$M_{+}=\begin{pmatrix}
1-x&-y&\\
-xy^{-1}&0
\end{pmatrix}$ 
and 
$M_{-}=\begin{pmatrix}
0&-x^{-1}y\\
-y^{-1}&1-x^{-1}\\
\end{pmatrix}$.
Let $|D|$ be the 4-valent graph obtained from $D$ by regarding each real crossing as a vertex of $|D|$. We 
denote by the same symbols $c_1, \dots, c_n$ the vertices of $|D|$. 
The graph $|D|$ is immersed in $\mathbb{R}^2$ and the multiple points of $|D|$ are virtual crossings of $D$. 
For each vertex $c_i$ of $|D|$, consider an open regular neighborhood $N(c_i, |D|)$ of $c_i$ in $|D|$. Then $N(c_i, |D|)-\{c_i\}$ is the union of four open arcs, which we call the {\it short edges around} $c_i$. According to the position, we denote by $i_0^-, i_1^-, i_0^+, i_1^+$ the short edges as in 
 Figure \ref{fig:labeledge}. 

We define a $2n\times 2n$ matrix, $P=(p_{kl})$ as follows. For each $i, j\in \{1, \dots, n\}$, 

$$p_{(2i-1+\epsilon)(2j-1+\lambda)}=\left\{
\begin{array}{ll}
1&\left(\parbox{10cm}{if two short edges  $i^{-}_{\epsilon}$ and $ j^{+}_{\lambda}$ are on the same edge of $|D|$}\right)\\
0&\parbox{3cm}{(otherwise)}
\end{array}\right. ,$$

 where  $\epsilon, \lambda\in \{0,1\}$.

\begin{thm}[Jaeger, Kauffman and Saleur \cite{rJKS}, Sawollek \cite{rSaw}]
For a virtual link diagram $D$, $Z_{D}(x, y) = (-1)^{w(D)}\mathrm{det}(M-P )$ is an invariant of the virtual link up to multiplication by powers of $x^{\pm 1}$, i.e., for any virtual link diagram $D'$ representing the same virtual link with $D$, we have $Z_{D'}(x, y)=x^mZ_D(x, y)$ for some $m\in \mathbb{Z}$.
\end{thm}

\begin{figure}[h]
\centerline{
\includegraphics[width=2.5cm]{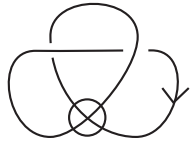}}
\caption{A virtual link diagram}\label{fig:exJKSSv1}
\end{figure}

The {\it JKSS invariant} of a virtual link $L$ is defined by $Z_L(x, y)=Z_D(x,y)$ for a diagram $D$ of $L$. 
For example, the JKSS invariant of a virtual link depicted in Figure~\ref{fig:exJKSSv1} is $(x-1)(y+1)(x+y)y^{-1}$.

We define an invariant of twisted links which is related to the JKSS invariant.
Let $D$ be a {twisted} link diagram with $n$ real crossings $c_1,\dots, c_n$.

We define a $4n\times 4n$  matrix $\widetilde{M}$, by $\widetilde{M}=\mathrm{diag}(\widetilde{M}_1,\dots, \widetilde{M}_n)$, where 
$\widetilde{M}_i=\widetilde{M}_+$ (or $\widetilde{M}_-$) if the crossing $c_i$ is positive (or negative). Here 
$\widetilde{M}_+$ and $\widetilde{M}_-$ are $4\times 4$ matrices $\mathrm{diag}(M_+, M_+)$ and $\mathrm{diag}(M_-, M_-)$ respectively, i.e., 
$$\widetilde{M}_+=
\begin{pmatrix}
1-x&-y&0&0\\
-xy^{-1}&0&0&0\\
0&0&1-x&-y\\
0&0&-xy^{-1}&0
\end{pmatrix}
\text{ and }
\widetilde{M}_-=
\begin{pmatrix}
0&-x^{-1}y&0&0\\
-y^{-1}&1-x^{-1}&0&0\\
0&0&0&-x^{-1}y\\
0&0&-y^{-1}&1-x^{-1}
\end{pmatrix}
.$$

For a twisted link diagram $D$, the graph $|D|$ is defined by the same way before. Each edge of $|D|$ may have bars on it. For each vertex $c_i$ of $|D|$, we denote by $i_0^-, i_1^-, i_0^+$, and $i_1^+$ the short edge around $c_i$ as before. 
We denote 
$i_{\epsilon}  \overset{\mathrm{e}}{\leftarrow} j_{\lambda}$ (or $i_{\epsilon} \overset{\mathrm{o}}{\leftarrow} j_{\lambda}$), 
if two short edges $i^{-}_{\epsilon}$ and $ j^{+}_{\lambda}$ are on the same edge of $|D|$ and there are an even (or odd) number of bars on the edge. 

We defined $4n\times 4n$ matrix, $\widetilde{P}=(\tilde{p}_{kl})$ as follows. For each $i, j\in\{1,\dots, n\}$, 
$$\tilde{p}_{(4i-3+a)(4j-3+b)}=\left\{
\begin{array}{ll}
1&\left(\parbox{9cm}{$i_a \overset{\mathrm{e}}{\leftarrow} j_b$ or $i_{3-a} \overset{\mathrm{e}}{\leftarrow} j_{3-b}$ or $i_a \overset{\mathrm{o}}{\leftarrow} j_{3-b}$ or $i_{3-a} \overset{\mathrm{o}}{\leftarrow} j_{b}$}\right)\\
0&\parbox{3cm}{(otherwise)}
\end{array}\right. , $$
where $a, b \in \{0, 1, 2, 3 \}$
Note that $i_k^-$ and $j_k^-$ arc not defined for $k\in\{2, 3\}$. We assume that $i_k \overset{\mathrm{e}}{\leftarrow} j_l$ and $i_k \overset{\mathrm{o}}{\leftarrow} j_l$ are false when  $k\in\{2, 3\}$ or  $l\in\{2, 3\}$.

\begin{thm}\label{thm1}
For a twisted link diagram $D$, $\widetilde{Z}_{D}(x, y) = \mbox{det}(\widetilde{M}-\widetilde{P })$ is an invariant of the twisted link up to multiplication by powers of $x^{\pm 1}$, i.e., for any twisted link diagram $D'$ representing the same twisted link with $D$, we have $Z_{D'}(x, y)=x^mZ_D(x, y)$ for some $m\in \mathbb{Z}$.
\end{thm}
For a twisted link $L$, we define the {\it twisted JKSS invariant} of $L$, denoted by $\widetilde{Z}_L(x, y)$, by $\widetilde{Z}_D(x, y)$ for a diagram $D$ of $L$.

\begin{figure}[!h]
\begin{center}
\begin{tabular}{cc}
\includegraphics[width=2cm]{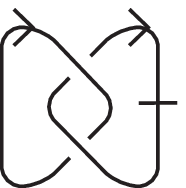}&
\includegraphics[width=2cm]{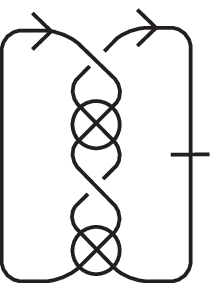}\\
(a)&(b)
\end{tabular}
\caption{The twisted JKSS invariants of twisted links}\label{fig:exJKSScalt1}
\end{center}
\end{figure}

The twisted JKSS invariant of the diagram in Figure \ref{fig:exJKSScalt1} (a) 
is $0$ and that of the diagram in Figure \ref{fig:exJKSScalt1} (b) is $y^{-2}\left(x^2-1\right) \left(y^2-1\right) \left(x^2-y^2\right)$. We see that they are not equivalent. 
Note that they are not distinguished by the twisted Jones polynomial defined in \cite{rBor}.  The twisted Jones polynomials of the twisted links  in Figure \ref{fig:exJKSScalt1} are  $-A^{-6}(A^{4}+A^{-4})M$. 

\section{Proof of Theorem \ref{thm1}}

Let $D$ be a twisted link diagram with bars $b_1,\dots , b_k$. 
Assume that $D$ is on the left of the $y$-axis and all bars are parallel to the $x$-axis with disjoint $y$-coordinates. Let $s(D)$ be  the twisted link diagram which is obtained from $D$ by reflection with respect to the $y$-axis and switching all real crossings of $D$.
See Figure~\ref{fig:doublecv1}.

We construct the double covering of $D$ as follows:

\begin{figure}
\centering{
\includegraphics[width=7cm]{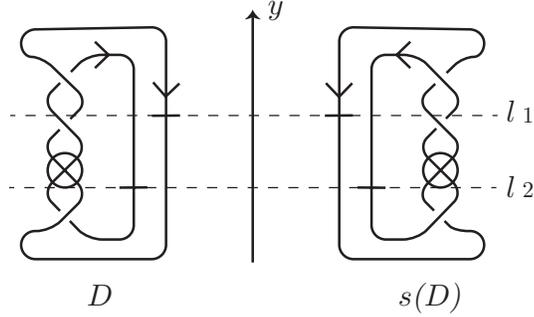}}
\caption{A twisted link diagram $D$ and $s(D)$}\label{fig:doublecv1}
\end{figure}

\begin{figure}[h]
\centering{
\includegraphics[width=10cm]{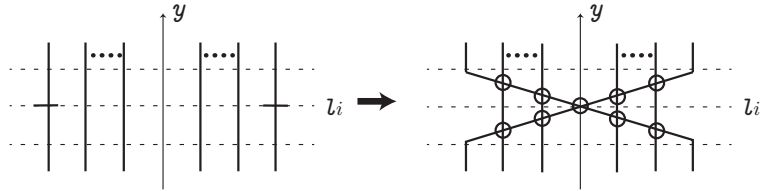} }
\caption{Construct double covering }\label{fig:doublecv2}
\end{figure}

For horizontal lines $l _1, \dots , l_k$ such that  $l_i$ contains $b_i$ and $s(b_i)$, 
we  replace  each part of $D\amalg s(D)$ in a neighborhood of $N(l_i)$  as in Figure~\ref{fig:doublecv2}. We call  this diagram $\widetilde{D}$ the {\it double covering diagram} of $D$.

\begin{figure}[h]
\centering{
\includegraphics[width=7cm]{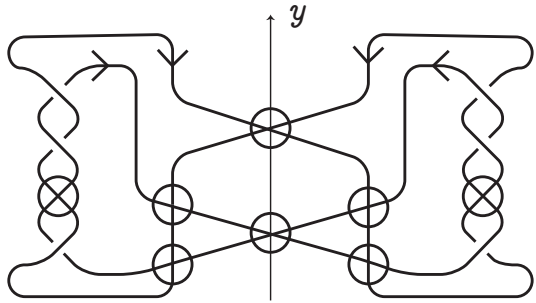} }
\caption{The double covering diagram of $D$}\label{fig:doublecv3}
\end{figure}

The double covering diagram of the twisted link diagram $D$ in Figure~\ref{fig:doublecv1} is shown in Figure~\ref{fig:doublecv3}. 

\begin{thm}[\cite{rkk7}]\label{thm2}
Let $D_1$ and $D_2$ be twisted link diagrams and 
$\widetilde{D_1}$ and $\widetilde{D_2}$ their double coverings diagrams of $D_1$ and $D_2$.
If $D_1$ and $D_2$ are equivalent as twisted links, then $\widetilde{D_1}$ and $\widetilde{D_2}$ are equivalent 
as virtual links.
\end{thm}

\begin{thm}\label{thm0}
For a twisted link diagram $D$, $\widetilde{Z}_{D}(x,y)$ coincides to ${Z}_{\widetilde{D}}(x,y)$, where $\widetilde{D}$ is the double covering of $D$.
\end{thm}

\noindent
{\bf Proof}\quad
Let $D$ be a twisted link diagram with $n$ real crossings $c_1, \cdots , c_n$. 
Let $s(c_1), \dots, s(c_n)$ be the real crossings of $s(D)$ corresponding to $c_1, \dots, c_n$. 
We regard $c_1, \cdots , c_n, s(c_1), \dots, s(c_n)$ as the real crossings of the double covering diagram $\widetilde{D}$. 
We rename the real crossings of $\widetilde{D}$ by $\tilde{c}_1, \cdots , \tilde{c}_{2n}$ such that $c_i=\tilde{c}_{2i-1}$ and $s(c_i)=\tilde{c}_{2i}$ for $i\in \{1,\dots, n\}$. 
Then the matrix $\widetilde{M}$ for $D$ coincides with matrix $M$ for $\widetilde{D}$. We show that the matrix $\widetilde{P}$ for $D$ coincides with the matrix $P$ for $\widetilde{D}$.
Suppose that real crossings $c_i$ and  $c_j$ of $D$ are the boundary points of an edge of $|D|$, i.e., $i_{\epsilon} \overset{\mathrm{e}}{\leftarrow} j_{\lambda}$ or  $i_{\epsilon} \overset{\mathrm{o}}{\leftarrow} j_{\lambda}$ holds for $D$. 

\begin{enumerate}
\item
If $i_{\epsilon} \overset{\mathrm{e}}{\leftarrow} j_{\lambda}$ for $D$, 
we see that ${c_i}$ and ${c_j}$ (or  $s(c_i)$  and ${s(c_j)}$ ) are the boundary points  of an edge of  $|\widetilde{D}|$. 
The short edges labeled by $(2i-1)_{\epsilon}^-$ and $(2j-1)_{\lambda}^+$ (or $(2i)_{1-\epsilon}^-$ and $(2j)_{1-\lambda}^+$) are on the same edge of $|\widetilde{D}|$. Thus we have 
$$
\tilde{p}_{(4i-3+a)(4j-3+b)}=
\left\{\begin{array}{ll}
1&(a=\epsilon \text{ and } b=\lambda)\\
1&(a=2+1-\epsilon \text{ and } b=2+1-\lambda)\\
0&(\text{otherwise}).
\end{array}\right.
$$ 
\item
If $i_{\epsilon} \overset{\mathrm{o}}{\leftarrow} j_{\lambda}$ for $D$, 
we see that ${c_i}$ and ${s(c_j)}$ (or  ${s(c_i)}$  and ${c_j}$ ) are the boundary points of an edge of $|\widetilde{D}|$. 
Namely, the short edges labeled by $(2i-1)_{\epsilon}^-$ and $(2j)_{1-\lambda}^+$ (or $(2i)_{1-\epsilon}^-$ and $(2j-1)_{\lambda}^+$) are on the same edge of $|\widetilde{D}|$. Thus we have 
$$
\tilde{p}_{(4i-3+a)(4j-3+b)}=
\left\{\begin{array}{ll}
1&(a=\epsilon \text{ and } b=2+1-\lambda)\\
1&(a=2+1-\epsilon \text{ and } b=\lambda)\\
0&(\text{otherwise}).
\end{array}\right.
$$ 
\end{enumerate}
Thus  we have the $\widetilde{P}=P$. Since $w(\widetilde{D})$ is $2w(D)$, $Z_{\widetilde{D}}(x, y)=\mathrm{det}(M-P)=\mathrm{det}(\widetilde{M}-\widetilde{P})
=\widetilde{Z}_D(x, y)$. 
\qed

%%%%%%%%%%%%%%%%%%%%%%%%%%%%%%%%%%%%%
\section{Properties of double covering diagrams}
For a twisted (virtual or classical) link diagram $D$, we denote by $r(D)$ the number of real crossings of $D$. 
For a twisted link (or a virtual link) $L$, we denote by $r(L)$ (or $r_v(L)$) the minimal number of real crossings of all of diagrams of $L$. 

\begin{thm}\label{prop1}
Let $L$ be a twisted link presented by a twisted link diagram $D$ and let $\widetilde{L}$ be the virtual link presented by the  double covering diagram $\widetilde{D}$ of $D$.
If $r(\widetilde{D})=r_v(\widetilde{L})$, then $r(D)=r(L)$.
\end{thm}
\noindent
{\bf Proof}\quad
Note that $r(\widetilde{D})=2r(D)$. 
If $r(D)\ne r(L)$, then there is a diagram of $L$, $D_0$ such that $r(D_0)=r(L)<r(D)$. 
For the double covering diagram $\widetilde{D}_0$ of $D_0$, we have $r(\widetilde{D}_0)=2r(D_0)<2r(D)$, 
which is a contradiction of the assumption since $\widetilde{D}_0$ is equivalent to $\widetilde{D}$ as a virtual link.  
\qed

\begin{prop}\label{prop2}
Let $D$ be  a twisted link diagram obtained from a virtual link diagram ${D_0}$ by adding a bar on an arc. Then $\widetilde{D}$ is equivalent to a connected sum of ${D_0}$ and $s({D_0})$ modulo virtual Reidemeister moves.
\end{prop}

\begin{figure}[h]
\centering{
\includegraphics[width=10cm]{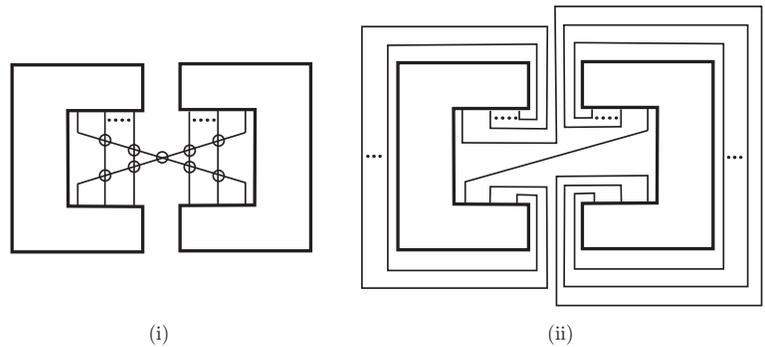}}
\caption{Double covering of a twisted link diagram with a bar}\label{fgprop2}
\end{figure}

\noindent
{\bf Proof}\quad
For the double covering $\widetilde{D}$ of $D$, there are two kinds of sets of virtual crossings. One is the set of virtual crossings which correspond to virtual crossings of  $D$ and $s(D)$, and the other 
is the set of virtual crossings which occur when we construct $\widetilde{D}$ from $D\amalg s(D)$ by replacement as in Figure~\ref{fig:doublecv2}. The second set of virtual crossings of $\widetilde{D}$ look as in Figure~\ref{fgprop2} (i).
As in Figure~\ref{fgprop2} (ii), such virtual crossings are eliminated by some virtual Reidemeister moves. 
The virtual link diagram shown in Figure~\ref{fgprop2} (ii) is a connected sum of two virtual link diagrams which are equivalent to ${D_0}$ and $s({D_0})$ modulo virtual Reidemeister moves.
\qed

\vspace{8pt}
%\vspace{3zw}
{\bf Acknowledgements }

The author  would like to thank Seiichi Kamada for his useful suggestion.

 \end{document}